\documentclass{amsart}
\usepackage{amssymb,latexsym}
\input epsf
\newtheorem{theorem}{Theorem}[section]

\newtheorem{corollary}[theorem]{Corollary}
\newtheorem{definition}[theorem]{Definition}

\newtheorem{lemma}[theorem]{Lemma}
\newtheorem{example}[theorem]{Example}

\newcommand{\rk}{{\rm rank}}
\newcommand{\nc}{{\rm NC}}

\newcommand{\RR}{{\mathbb R}}
\newcommand{\ZZ}{{\mathbb Z}}

\newcommand{\bB}{{\mathcal B}}
\newcommand{\eE}{{\mathcal E}}
\newcommand{\lleft}{( \! (}
\newcommand{\rright}{) \! )}
\renewcommand{\to}{\rightarrow}

\begin{document}
\title[Shellability of noncrossing partition lattices]
{Shellability of noncrossing partition lattices}

\author{Christos~A.~Athanasiadis}
\address{Department of Mathematics\\
University of Crete\\
71409 Heraklion, Crete, Greece}
\email{caa@math.uoc.gr}

\author{Thomas~Brady}
\address{School of Mathematical Sciences\\
Dublin City University\\
Glasnevin, Dublin 9\\
Ireland}
\email{tom.brady@dcu.ie}

\author{Colum~Watt}
\address{School of Mathematics\\
Dublin Institute of Technology\\
Dublin 8\\
Ireland}
\email{colum.watt@dit.ie}

\date{February 28, 2005}
\thanks{
Supported in part by the American Institute of Mathematics (AIM) 
and the NSF}
\begin{abstract}
We give a case-free proof that the lattice of noncrossing partitions
associated to any finite real reflection group is EL-shellable. 
Shellability of these lattices was open for the groups of type $D_n$ 
and those of exceptional type and rank at least three.
\end{abstract}

\maketitle

\section{Introduction}
\label{intro}

Consider a finite real reflection group $W$ and the partial order on $W$ 
defined by letting $u \preceq v$ if there exists a shortest factorization 
of $u$ as a product of reflections in $W$ which is a prefix of such a 
shortest factorization of $v$. This order turns $W$ into a graded poset 
having the identity $1$ as its unique minimal element, where the rank of 
$w$ is the length of a shortest factorization of $w$ into reflections. 
For any $w \in W$ we denote by $\nc_W (w)$ the interval $[1, w]$ in 
this partial order. We are primarily interested in the case that $w$ 
is a Coxeter element $\gamma$ of $W$, viewed as a finite Coxeter group. 
Since all Coxeter elements of $W$ are conjugate to each other, the 
isomorphism type of the poset $\nc_W (\gamma)$ is independent of $\gamma$. 
We denote this poset by $\nc_W$ when the choice of the Coxeter element 
$\gamma$ is irrelevant and call $\nc_W$ the \emph{noncrossing partition
lattice} associated to $W$. 

The poset $\nc_W$ plays a crucial role in the construction of new monoid 
structures and $K(\pi, 1)$ spaces for Artin groups associated with finite 
Coxeter groups; see for instance \cite{Be, Br, BW1} and the survey articles 
\cite{Mc1, Mc2}. It is self-dual \cite[Section 2.3]{Be}, graded and has 
been verified case by case to be a lattice \cite[Fact 2.3.1]{Be} 
\cite[Section 4]{BW1}. For Coxeter groups of types $A$ and $B$ in the 
Cartan-Killing classification it is isomorphic to the classical lattice of 
noncrossing partitions of an $n$-element set, defined and studied by 
Kreweras \cite{Kr}, and to its type $B$ analogue, defined by Reiner 
\cite{Re}, respectively. The main purpose of the present article is to 
give a case-free proof of the following theorem.

\begin{theorem} 
The poset $\nc_W$ is EL-shellable for any finite Coxeter group $W$. 
\label{thm:main}
\end{theorem}

EL-shellability (see Section \ref{pre}) for a bounded graded poset $P$ of 
rank $r$ implies that the simplicial complex $\Delta(\bar{P})$ of chains in 
the proper part of $P$ is shellable. In particular it implies that the 
geometric realization of $\Delta(\bar{P})$ has the homotopy type of a wedge 
of $(r-2)$-dimensional spheres, the number of which is determined by the 
EL-labeling, and that the Stanley-Reisner ring of $\Delta(\bar{P})$ is 
Cohen-Macaulay over an arbitrary field. 

Theorem \ref{thm:main} was proved by Bj\"orner and Edelman \cite[Example 
2.9]{Bj} in the case of Coxeter groups of type $A$ and 
by Reiner \cite[Section 6]{Re} in the case of type $B$ while it was left 
open in \cite[Section 7]{AR} in that of type $D$ and was conjectured by 
Reiner \cite{Re02} for all finite Coxeter groups. After introducing basic 
definitions, notation and some preliminary lemmas related to real 
reflection groups, root systems and shellability (Section \ref{pre}) we 
proceed as follows. We give sufficient conditions on a total ordering of 
the set of reflections $T$ and a Coxeter element $\gamma$ of a 
crystallographic group $W$ for the natural edge labeling of $\nc_W 
(\gamma)$ with label set $T$ to be an EL-labeling (Theorem \ref{thm:el}). 
We exploit these conditions to describe explicit families of EL-labelings 
in the cases of the classical reflection groups (Examples \ref{ex:A} 
and \ref{ex:BD} and Corollary \ref{cor:ABD}). The existence of a total 
ordering and Coxeter element for any $W$ which satisfy these sufficient 
conditions follows from recent work of the last two authors \cite{BW2}, 
where a case-free proof of the lattice property of $\nc_W$ is given. The 
particular orderings and Coxeter elements considered there (see Section 
\ref{stein}) were introduced by Steinberg \cite{Ste} and play a 
crucial role in the constructions of \cite{BW2}. It is shown in Section
\ref{stein} (Theorem \ref{thm:any}) that they provide EL-labelings of 
$\nc_W$ for any finite real
reflection group (including the non-crystallographic ones). A result on 
the M\"obius function of $\nc_W$ is deduced (Corollary \ref{cor:mobius}). 

\section{Preliminaries}
\label{pre} 
 
\noindent
{\bf Reflection groups and root systems.}
Let $W$ be a finite real reflection group of rank $\ell$ with a corresponding 
root system $\Phi$, acting faithfully by orthogonal transformations on an 
$\ell$-dimensional Euclidean space $V$ with inner product $( \ , \, )$. 
Thus $W$ is generated by elements acting as orthogonal reflections in $V$ 
while $\Phi$ is invariant under the action of $W$ and consists of a pair 
$\{\alpha, -\alpha\}$ of nonzero vectors for each reflection $t$ in $W$, 
which are orthogonal to the reflecting hyperplane of $t$. We denote by $T$ 
the set of all reflections in $W$, by $H_\alpha$ the linear hyperplane in $V$
orthogonal to $\alpha$ and by $t_\alpha$ the orthogonal reflection in 
$H_\alpha$. We refer the reader to the texts by Bj\"orner 
and Brenti \cite{BB} and Humphreys \cite{Hu} for any undefined terminology 
and background on reflection groups and root systems. In particular we assume 
the notions of a positive system, simple system and simple reflection. A 
\emph{Coxeter element} of $W$ is the product of the simple reflections in an 
arbitrary order for any choice of simple system for $\Phi$. The order of any 
Coxeter element of $W$ is \emph{the Coxeter number}, denoted by $h$. The 
\emph{root lattice} $Q_\Phi$ is the $\ZZ$-span of $\Phi$. The root system 
$\Phi$ is called \emph{crystallographic} and $W$ is called a \emph{Weyl 
group} if $W$ preserves $Q_\Phi$. 

The next lemma follows from the results of Sommers in \cite[Section 3]{So} 
(see also \cite[Lemma 2.3]{Pa}).
\begin{lemma} {\rm (\cite{So})}
If $\alpha_1, \alpha_2,\dots,\alpha_r \in \Phi$ and $\alpha_1 + \alpha_2 
+ \cdots + \alpha_r = \alpha \in \Phi$ then $\alpha_1 = \alpha$ or there 
exists $i$ with $2 \le i \le r$ such that $\alpha_1 + \alpha_i \in \Phi 
\cup \{0\}$.
\label{lem:ro}
\end{lemma}
\begin{proof}
See \cite[Lemma 2.1 (ii)]{Ath}.
\end{proof}
\begin{lemma}
If $\{\alpha_1, \alpha_2,\dots,\alpha_\ell\} \subseteq \Phi^+$ is a 
$\ZZ$-basis of the root lattice $Q_\Phi$ and $\alpha_i - \alpha_j \notin 
\Phi$ for all $i < j$ then $\{\alpha_1, \alpha_2,\dots,\alpha_\ell\}$ is 
the set of simple roots in $\Phi^+$.
\label{lem:simple}
\end{lemma}
\begin{proof}
It suffices to prove that any positive root $\alpha \in \Phi^+$ can be 
written as a linear combination of $\alpha_1, 
\alpha_2,\dots,\alpha_\ell$ with nonnegative integer coefficients. Since 
$\{\alpha_1, \alpha_2,\dots,\alpha_\ell\}$ is a $\ZZ$-basis of the root 
lattice we can write uniquely $\alpha = r_1 \alpha_1 + r_2 \alpha_2 + 
\cdots + r_\ell \alpha_\ell$ for some integers $r_i$. Let us rewrite 
this expression as
\[ \alpha = \beta_1 + \beta_2 + \cdots + \beta_r \]
where there are $r_i$ copies of $\alpha_i$ or $-r_i$ copies of $-\alpha_i$ 
among the $\beta$'s if $r_i \ge 0$ or $r_i < 0$, respecively. To show that 
$r_i \ge 0$ for $1 \le i \le \ell$ we must show that $\beta_i \in \Phi^+$ 
for $1 \le i \le r$. Suppose that $\beta_i$ is a negative root for some 
$i$, say for $i=1$. Since $\alpha \in \Phi^+$ at least one of the $\beta$'s 
is a positive root. By applying repeatedly Lemma \ref{lem:ro} and 
reordering the $\beta$'s if necessary we may assume that there exists an 
index $2 \le j \le r$ such that $\beta_1 + \beta_2 + \cdots + \beta_i \in 
\Phi$ for $1 \le i < j$, $\beta_j \in \Phi^+$, $\beta_i \in -\Phi^+$ for 
$1 \le i < j$ and $\beta_1 + \beta_2 + \cdots + \beta_j \in \Phi \, \cup 
\{0\}$. Since the $\alpha$'s are linearly independent we must have 
$\beta_1 + \beta_2 + \cdots + \beta_j \in \Phi$. Applying Lemma 
\ref{lem:ro} once more we conclude that $\beta_i + \beta_j \in \Phi$ for 
some $1 \le i < j$, which means that $\alpha_i - \alpha_j \in \Phi$ for 
some $i<j$, contrary to the hypothesis.
\end{proof}

\medskip
\noindent
{\bf EL-labelings and shellability.} Let $(P, \le)$ be a finite bounded 
graded poset (short for partially ordered set). Thus $P$ has a unique minimal 
and a unique maximal element, denoted $\hat{0}$ and $\hat{1}$ respectively, 
and all maximal (with respect to inclusion) chains in $P$ have the same 
length (one less than their number of elements), called the \emph{rank} 
of $P$ and denoted $\rk(P)$. Let $\eE(P)$ 
be the set of covering relations of $P$, meaning pairs $(x, y)$ of elements 
of $P$ such that $x < y$ holds in $P$ and $x < z \le y$ holds only for 
$z = y$, and let $\Lambda$ be a totally ordered set. An \emph{edge 
labeling} of $P$ with label set $\Lambda$ is a map $\lambda: \eE(P) \to 
\Lambda$. Let $c$ be an unrefinable chain $x_0 < x_1 < \cdots < x_r$ of 
elements of $P$, so that $(x_{i-1}, x_i) \in \eE(P)$ for all $1 \le i 
\le r$. We let $\lambda(c) = (\lambda(x_0, x_1), \lambda(x_1, 
x_2),\dots,\lambda(x_{r-1}, x_r))$ be the label of $c$ with respect to 
$\lambda$ and call $c$ \emph{rising} or \emph{falling} with respect to 
$\lambda$ if the entries of $\lambda(c)$ strictly increase or weakly 
decrease, respectively, in the total order of $\Lambda$. We say that $c$
is \emph{lexicographically smaller} than an unrefinable chain $c'$ in $P$ 
of the same length (with respect to $\lambda$) if $\lambda(c)$ precedes 
$\lambda(c')$ in the lexicographic order induced by the total order of 
$\Lambda$. 
\begin{definition} {\rm (\cite{Bj})}
An edge labeling $\lambda$ of $P$ is called an EL-labeling if for every 
non-singleton interval $[u, v]$ in $P$
\begin{enumerate}
\itemsep=0pt
\item[(i)] there is a unique rising maximal chain in $[u, v]$ and
\item[(ii)] this chain is lexicographically smallest among all maximal 
chains in $[u, v]$
\end{enumerate}
with respect to $\lambda$. The poset $P$ is called EL-shellable if it has 
an EL-labeling for some label set $\Lambda$.  
\label{def:el}
\end{definition}
See \cite[Chapter 3]{Sta} for more background on partially ordered sets 
and EL-shellability. If $P$ is EL-shellable then the simplicial complex 
of chains (order complex) $\Delta (\bar{P})$ of the proper part $\bar{P} 
= P - \{\hat{0}, \hat{1}\}$ of $P$ has the homotopy type of a wedge of 
spheres of dimension $\rk(P) - 2$. The number of these spheres is equal
to the number of falling maximal chains of $P$ with respect to the 
EL-labeling and is also equal to the M\"obius number of $P$, up to the 
sign $(-1)^{\rk(P)}$; see \cite[Section 3.13]{Sta}. 

\medskip
\noindent
{\bf Reflection length and the poset $\nc_W (w)$.} For $w \in W$ 
let $l_T (w)$ denote the smallest $k$ such that $w$ can be written as a  
product of $k$ reflections in $T$. The partial order $\preceq$ on $W$ is 
defined by letting 
\[ u \preceq v \ \ \ \text{if and only if} \ \ \ l_T (u) + l_T (u^{-1} 
v) = l_T (v), \]
in other words if there exists a shortest factorization of $u$ into 
reflections in $T$ which is a prefix of such a shortest factorization of 
$v$. With this order $W$ is a graded poset having the identity $1$ as 
its unique minimal element and rank function $l_T$. We denote by $F(w)$
the subspace of $V$ fixed by $w$. The next lemma follows from the results 
in \cite[Section 2]{BW1} (see also \cite[Lemma 1.2.1]{Be}). 
\begin{lemma} {\rm (\cite{BW1})}
Let $w \in W$. 
\begin{enumerate}
\itemsep=0pt
\item[(i)] $l_T (w) = \ell - \dim F(w)$. 
\item[(ii)] For $\alpha \in \Phi^+$ we have $t_\alpha \preceq w$ if and 
only if $F(w) \subseteq H_\alpha$. 
\end{enumerate}
\label{lem0}
\end{lemma}
The interval $[1, w]$ in $(W, \preceq)$, denoted $\nc_W (w)$, is also graded 
with rank function $l_T$ and has rank $\ell = \dim V$ if $w$ is a Coxeter 
element of $W$. 

\section{Reflection Orderings, Coxeter elements and EL-Labelings}
\label{orderings}

Let $W$ be a finite real reflection group of rank $\ell$ with corresponding 
root system $\Phi \subset V$ and set of reflections $T$. Let $\Phi^+ 
\subset \Phi$ be a fixed choice of a positive system. A total ordering $<$ 
of $T$ is called a \emph{reflection ordering} for $W$ \cite[Section 
5.1]{BB} if whenever $\alpha, \alpha_1, \alpha_2 \in \Phi^+$ are distinct
roots and $\alpha$ is a positive linear combination of $\alpha_1$ and 
$\alpha_2$ we have either \[ t_{\alpha_1} < t_\alpha < t_{\alpha_2} \] or 
\[ t_{\alpha_2} < t_\alpha < t_{\alpha_1}. \]       

\smallskip
An induced subsystem of $\Phi$ of rank $i$ is the intersection $\Phi'$ of 
$\Phi$ with the linear span of $i$ linearly independent roots in $\Phi$. 
The set $\Phi'$ is a root system on its own and $\Phi' \cap  \Phi^+$ is 
a choice of a positive system for $\Phi'$. The following is the  main 
definition in this section.
\begin{definition}
A reflection ordering $<$ of $T$ is compatible with a Coxeter element 
$\gamma$ of $W$ if for any irreducible rank 2 induced subsystem $\Phi' 
\subseteq \Phi$ the following holds: if $\alpha$ and $\beta$ 
are the simple roots of $\Phi'$ with respect to $\Phi' \cap \Phi^+$ and 
$t_\alpha t_\beta \in \nc_W (\gamma)$ then $t_\alpha < t_\beta$. 
\label{def:comp}
\end{definition}
Observe that a rank 2 subsystem $\Phi' \subseteq \Phi$ is irreducible 
if and only if $\Phi' \cap \Phi^+$ has at least three roots.
\begin{example} {\rm
Let $W$ be of rank 2, so that $W$ is a dihedral group of order $2m$ for 
some $m \ge 2$. Let $\alpha$ and $\beta$ be the two simple roots in $\Phi$. 
There are $m$ reflections in $W$, namely $t_1, t_2,\dots,t_m$ where $t_i 
= t_\alpha (t_\beta t_\alpha)^{i-1}$ (so $t_1 = t_\alpha$ and $t_m 
= t_\beta$) and two reflection orderings, namely $t_1 < t_2 < \cdots 
< t_m$ and its reverse. The reflection ordering $<$ is compatible 
with the Coxeter element $\gamma = t_\alpha t_\beta$. The poset $\nc_W 
(\gamma)$ has $m$ maximal chains corresponding to the $m$ shortest 
factorizations $\gamma = t_1 t_m = t_2 t_1 = \cdots = t_m t_{m-1}$ of 
$\gamma$ into reflections. Clearly the chain corresponding to the 
factorization $\gamma = t_1 t_m$ is the unique rising maximal chain and 
the lexicographically smallest maximal chain in the edge labeling of 
$\nc_W (\gamma)$ which assigns the label $(t_i, t_j)$ to the chain 
corresponding to a factorization $\gamma = t_i t_j$, when the label set 
$T$ is totally ordered by $<$. \qed
}
\label{ex:dihedral}
\end{example}
\begin{example} {\rm
Let $W$ be of type $A_{n-1}$ and $\Phi^+ = \{e_i - e_j: 1 \le i < j \le 
n\}$ be the set of positive roots. Thus $W$ is isomorphic to the symmetric 
group of permutations of the set $\{1, 2,\dots,n\}$ and the reflections 
correspond to the transpositions $(i, j)$ for $1 \le i < j \le n$. If the
$n$-cycle $\gamma = (1, 2,\dots,n)$ is chosen as the Coxeter element for $W$
then there is a canonical isomorphism of $\nc_W (\gamma)$ with the classical
lattice of noncrossing partitions of $\{1, 2,\dots,n\}$ (see for instance 
\cite[Section 4]{Be} or \cite[Section 3]{Br}). One can easily check that 
any total ordering on the set of transpositions $(i, j)$, where $1 
\le i < j \le n$, for which $$(i, j) < (i, k) < (j, k)$$ for $1 \le i < j < 
k \le n$ (such as the lexicographic ordering) induces a reflection ordering 
for $W$ which is compatible with $\gamma$. 
\qed
}
\label{ex:A}
\end{example}
\begin{example} {\rm
Let $W$ be of type $B_n$ or $D_n$ and 
\[ \Phi^+ = \begin{cases}
    \{e_i: 1 \le i \le n\} \, \cup \, \{e_i \pm e_j: 1 \le i < j \le n\}, 
    & \text{for $W$ of type $B_n$}\\
    \{e_i \pm e_j: 1 \le i < j \le n\}, 
    & \text{for $W$ of type $D_n$}   
\end{cases} \]
be the set of positive roots. In what follows we use the notation of 
\cite[Section 3]{BW1} (and \cite[Section 2]{AR}) to represent elements of 
$W$ as certain permutations of the set $\{1, 2,\dots,n, -1, -2,\dots,-n\}$, 
so that $\lleft i_1, i_2,\dots,i_k \rright$ stands for the product of cycles 
$(i_1, i_2,\dots,i_k)(-i_1, -i_2,\dots,-i_k)$ and $[i_1, i_2,\dots,i_k]$ 
for the ballanced cycle $(i_1, i_2,\dots,i_k, -i_1, -i_2,\dots,-i_k)$. If 
\[ \gamma = \begin{cases}
    [1, 2,\dots,n], & \text{for $W$ of type $B_n$}\\
    [1, 2,\dots,n-1] \, [n], & \text{for $W$ of type $D_n$}   
\end{cases} \]
is chosen as the Coxeter element of $W$ then there is a canonical 
isomorphism of $\nc_W (\gamma)$ with Reiner's $B_n$ analogue \cite{Re} of 
the lattice of noncrossing partitions of $\{1, 2,\dots,n\}$ (see 
\cite[Section 4]{Be} or \cite[Section 3]{BW1}) if $W$ is of type $B_n$, 
and an explicit combinatorial description of $\nc_W (\gamma)$ in terms of 
planar noncrossing diagrams \cite[Section 3]{AR} if $W$ is of type $D_n$.
One can check directly that in the case of type $B_n$ any total ordering on 
the set of signed transpositions $\lleft i, j \rright$ and $\lleft i,-j 
\rright$ for $1 \le i < j \le n$ and $[i] = (i,-i)$ for $1 \le i \le n$ 
for which
\[ \begin{tabular}{ll}
$\lleft i, j \rright < \lleft i, k \rright < \lleft j, k \rright$ & for $1 \le 
i < j < k \le n$ \\ 
$\lleft i, j \rright < \lleft i, -k \rright < \lleft j, -k \rright$ & for $1 
\le i < j \le n$ and $1 \le k < i$ or $j < k \le n$ \\ 
$\lleft i, j \rright < [i] < \lleft i, -j \rright < [j]$ & for $1 \le i < j 
\le n$, 
\end{tabular} \]
such as the ordering
\[ \begin{tabular}{r}
$\lleft 1,2 \rright < \lleft 1,3 \rright < \cdots < \lleft 1,n \rright < 
\lleft 2,3 \rright < \cdots < \lleft n-1,n \rright <$ \\ 
$[1] < \lleft 1,-2 \rright < \lleft 1,-3 \rright < \cdots < \lleft 1,-n 
\rright <$ \\ 
$< [2] < \lleft 2,-3 \rright < \cdots < \lleft 2,-n \rright <$ \\ 
$\cdots < [n-1] < \lleft n-1,-n \rright < [n]$,
\end{tabular} \]
induces a reflection ordering for $W$ which is compatible with $\gamma$. 
Similarly in the case of type $D_n$ any total ordering on the set of signed 
transpositions $\lleft i, j \rright$ and $\lleft i,-j \rright$ for $1 \le 
i < j \le n$ for which
\[ \begin{tabular}{ll}
$\lleft i, j \rright < \lleft i, k \rright < \lleft j, k \rright$ & for $1 
\le i < j < k \le n$ \\ 
$\lleft i, j \rright < \lleft i, -k \rright < \lleft j, -k \rright$ & for $1 
\le i < j \le n-1$ and $1 \le k < i$ or $j < k \le n$ \\ 
$\lleft i, \pm n \rright < \lleft i, -j \rright < \lleft j, \pm n \rright$ 
& for $1 \le i < j \le n-1$,  
\end{tabular} \]
such as the ordering
\[ \begin{tabular}{r}
$\lleft 1,2 \rright < \lleft 1,3 \rright < \cdots < \lleft 1,n-1 \rright < 
\lleft 2,3 \rright < \cdots < \lleft n-2,n-1 \rright <$ \\ 
$\lleft 1,n \rright < \lleft 1,-n \rright < \lleft 1,-2 \rright < \lleft 1,
-3 \rright < \cdots < \lleft 1,-(n-1) \rright <$ \\ 
$\lleft 2,n \rright < \lleft 2,-n \rright < \lleft 2,-3 \rright < \cdots < 
\lleft 2,-(n-1) \rright <$ \\ $\cdots < 
\lleft n-1,n \rright < \lleft n-1,-n \rright$,
\end{tabular} \]
induces a reflection ordering for $W$ which is compatible with $\gamma$. 
\qed
}
\label{ex:BD}
\end{example}

Let $\gamma$ be a Coxeter element of $W$. Any covering relation in $\nc_W 
(\gamma)$ is of the form $(u, v)$ with $u^{-1} v \in T$. Setting $\lambda 
(u, v) = u^{-1} v$ for any such covering relation defines an edge labeling 
$\lambda$ of $\nc_W (\gamma)$ with label set $T$ which we call the 
\emph{natural edge labeling} of $\nc_W (\gamma)$. Part (ii) of the 
following theorem is the main result of this section.
\begin{theorem}
Let $W$ be a finite real reflection group with set of reflections $T$ 
and Coxeter element $\gamma$ and let $\lambda$ be the natural edge 
labeling of $\nc_W (\gamma)$.
\begin{enumerate}
\itemsep=0pt
\item[(i)] For any total ordering of $T$ and any non-singleton interval 
$[u, v]$ in $\nc_W (\gamma)$ there is a unique lexicographically smallest 
maximal chain in $[u, v]$ and this chain is rising with respect to
$\lambda$.
\item[(ii)] If $T$ is totally ordered by a reflection ordering which 
is compatible with $\gamma$ and $W$ is a Weyl group then $\lambda$ is an 
EL-labeling. 
\end{enumerate} 
\label{thm:el}
\end{theorem}
We first need to establish a few lemmas. For any interval 
$[u, v]$ in $\nc_W (\gamma)$ we denote by $\lambda ([u,v])$ the set of 
natural labels of all maximal chains in $[u, v]$, with the convention that 
$\lambda([u, v]) = \{\emptyset\}$ if $u=v$, and abbreviate $\lambda 
([1,w])$ as $\lambda(w)$.
\begin{lemma}
Let $[u, v]$ be a non-singleton interval in $\nc_W (\gamma)$. 
\begin{enumerate}
\itemsep=0pt
\item[(i)] If $[u, v]$ has length two and $(s, t) \in \lambda ([u,v])$ then 
$(t, s') \in \lambda ([u,v])$ for some $s' \in T$.
\item[(ii)] It $t \in T$ appears in some coordinate of an element of $\lambda 
([u,v])$ then $t = \lambda(u, u')$ for some covering relation $(u, u')$ in 
$[u, v]$. 
\item[(iii)] The reflections appearing as the coordinates of an element of 
$\lambda ([u,v])$ are pairwise distinct. 
\end{enumerate}
\label{lem1}
\end{lemma}
\begin{proof}
Part (i) is clear since $T$ is closed under conjugation, so that $s' = tst$ 
is also an element of $T$. Parts (ii) and (iii) follow from repeated 
application of part (i).
\end{proof}

\begin{lemma} 
Let $[u, v]$ be a non-singleton interval in $\nc_W (\gamma)$ and let $w = 
u^{-1} v$. There is a poset isomorphism $f: [1, w] \to [u, v]$ such that
$\lambda (x,y) = \lambda (f(x), f(y))$ for all covering relations $(x, y)$ 
in $[1, w]$.
\label{lem4}
\end{lemma}
\begin{proof}
It follows immediately from the definitions that the map $f: [1, w] \to [u, 
v]$ with $f(x) = ux$ for $x \in [1, w]$ is well defined and has the desired 
properties.
\end{proof}

Recall that a \emph{parabolic Coxeter element} in $W$ is a Coxeter element 
$w$ in a parabolic subgroup of $W$, meaning a subgroup generated by a subset 
of a set of Coxeter generators of $W$. This subgroup, denoted $W_w$, depends 
only on $w$ and contains all reflections $t \preceq w$ (see \cite[Corollary 
1.6.2]{Be}).
\begin{lemma} {\rm (\cite[Lemma 1.4.3]{Be})}
Let $w \in W$. There exists a Coxeter element $\gamma$ of $W$ with $w 
\preceq \gamma$ if and only if $w$ is a parabolic Coxeter element.
\qed
\label{lem2}
\end{lemma}
\begin{lemma}
If $w \preceq \gamma$ for some Coxeter element $\gamma$ of $W$ then any 
reflection ordering for $W$ which is compatible with $\gamma$ restricts to 
a reflection ordering for $W_w$ which is compatible with $w$.
\label{lem5}
\end{lemma}
\begin{proof}
We may assume that $w$ has rank at least two in $\nc_W (\gamma)$. Let 
$\Phi_w \subseteq \Phi$ be the root system corresponding to $W_w$. We 
first check that if $\Phi'_w$ is an induced rank 2 subsystem of $\Phi_w$ 
and $\Phi' = \RR \Phi'_w \cap \Phi$ is the corresponding induced rank 2 
subsystem of $\Phi$ then $\Phi'_w = \Phi'$. Clearly $\Phi'_w \subseteq 
\Phi'$. Let $\alpha$ and $\beta$ be the two simple roots of $\Phi'_w$ 
and let $\alpha' \in \Phi'$. To show that $\alpha' \in \Phi'_w$ note 
that $H_\alpha \cap H_\beta \subset H_{\alpha'}$ and hence $F(w) \subset 
H_{\alpha'}$. It follows from Lemma \ref{lem0} (ii) that $t_{\alpha'} 
\preceq w$ and hence that $t_{\alpha'} \in W_w$, in other words that 
$\alpha' \in \Phi'_w$. To conclude the proof suppose that $\Phi'_w$ is 
irreducible and that $t_\alpha t_\beta \preceq w$. 
From $w \preceq \gamma$ we conclude that $t_\alpha t_\beta \preceq 
\gamma$. Since $\Phi'_w = \Phi'$ and the reflection ordering on $T$ 
is compatible with $\gamma$, $t_\alpha$ must preceed $t_\beta$ in this 
ordering. 
\end{proof}

Recall also from \cite[Section 1.5]{Be} that the Artin group $\bB_\ell$ of 
type $A_{\ell-1}$ acts on the set of shortest factorizations of $\gamma$ into 
reflections or, equivalently, on the set $\lambda(\gamma)$. More precisely,
the $i$th generator of $\bB_\ell$ acts on $(t_1, t_2,\dots,t_\ell) \in 
\lambda(\gamma)$ by replacing the pair $(t_i, t_{i+1})$ by $(t_i t_{i+1} 
t_i, t_i)$ while leaving other coordinates of $\lambda(\gamma)$ fixed.
\begin{lemma} 
The action of $\bB_\ell$ on the set of shortest factorizations of $\gamma$ 
into reflections is transitive. 
\label{lem3}
\end{lemma}
\begin{proof}
See Proposition 1.6.1 in \cite{Be}.
\end{proof}

\begin{lemma}
Let $t_{\alpha_1} t_{\alpha_2} \cdots t_{\alpha_\ell}$ be a shortest
factorization of a Coxeter element of $W$ into reflections.
\begin{enumerate}
\itemsep=0pt
\item[(i)] $\{\alpha_1, \alpha_2,\dots,\alpha_\ell\}$ is a linear
basis of $V$.
\item[(ii)] If $W$ is a Weyl group then $\{\alpha_1, 
\alpha_2,\dots,\alpha_\ell\}$ is a $\ZZ$-basis of the root lattice 
$Q_\Phi$.
\end{enumerate}
\label{lem:basis}
\end{lemma}
\begin{proof}
Part (i) is follows from the fact that Coxeter elements have trivial
fixed space in $V$.
The conclusion of part (ii) is clear for a shortest factorization of the 
Coxeter element into simple reflections. In view of Lemma \ref{lem3} it 
suffices to show that if two shortest factorizations $t_{\alpha_1} 
t_{\alpha_2} \cdots t_{\alpha_\ell}$ and $t_{\beta_1} t_{\beta_2} \cdots 
t_{\beta_\ell}$ are related by the action of a single generator of the 
Artin group $\bB_\ell$ then $\{\alpha_1, \alpha_2,\dots,\alpha_\ell\}$ is 
a $\ZZ$-basis of $Q_\Phi$ if and only if the same is true for 
$\{\beta_1, \beta_2,\dots,\beta_\ell\}$. We may thus assume that there 
exists an index $1 \le i < \ell$ such that $\alpha_j = \beta_j$ for $j 
\neq i, i+1$, $\beta_{i+1} = \alpha_i$ and $t_{\beta_i} = t_{\alpha_i} 
t_{\alpha_{i+1}} t_{\alpha_i}$. The last equality implies that
\[ \pm \beta_i = t_{\alpha_i} (\alpha_{i+1}) = \alpha_{i+1} - 
\frac{2 (\alpha_i, \alpha_{i+1})}{(\alpha_i, \alpha_i)} \, \alpha_i \]
and the claim follows since $2 (\alpha_i, \alpha_{i+1}) / (\alpha_i, 
\alpha_i) \in \ZZ$. 
\end{proof}

\vspace{0.1 in}
\noindent
\emph{Proof of Theorem \ref{thm:el}.} In what follows we will write 
$\nc_W$ instead of $\nc_W (\gamma)$. 

\medskip
\noindent
(i) We proceed by induction on the length of the interval $[u, v]$, 
the claim being trivial if this is equal to one. Clearly all covering 
relations of the form $(u, u')$ in $[u,v]$ have distinct labels. Let 
$(u, ut)$ be the one with the smallest label $t$. In view of the 
induction hypothesis applied to the interval $[ut, v]$ it suffices to 
prove that all covering relations in $[ut, v]$ have labels greater 
than $t$. This follows from parts (ii) and (iii) of Lemma \ref{lem1} 
which imply that any such label is different from $t$ and equal to the 
label of a covering relation $(u, u')$ in $[u,v]$.

\medskip
\noindent
(ii) Let $<_\gamma$ be the reflection 
ordering of $T$ which is compatible with the Coxeter element $\gamma$ 
of $W$. In view of part (i) it remains to show that there is at most one 
rising maximal chain in a non-singleton interval $[u, v]$ with respect 
to $\lambda$. In view of Lemma \ref{lem4} it suffices to prove this for 
an interval of the form $[1, w]$ in $\nc_W$. Since $w \preceq 
\gamma$, by Lemma \ref{lem2} $w$ is a parabolic Coxeter element in $W$. 
By Lemma \ref{lem5} the restriction of $<_\gamma$ on the set of 
reflections of the parabolic subgroup $W_w$ is a reflection ordering 
which is compatible with $w$. Hence we may as well assume that $w = \gamma$, 
so that $[1, w]$ is the entire poset $\nc_W = \nc_W (\gamma)$. Clearly there 
is at most one rising maximal chain in $\nc_W$ whose label is a permutation 
of the set $S$ of simple reflections. Hence it suffices to 
show that a maximal chain $c$ in $\nc_W$ whose label $\lambda(c)$ is not 
a permutation of $S$ cannot be rising. Indeed, let $c$ be such a maximal 
chain and let $\lambda (c) = (t_{\alpha_1}, 
t_{\alpha_2},\dots,t_{\alpha_\ell})$. By Lemmas \ref{lem:simple} and 
\ref{lem:basis} (ii) there exist indices $i < j$ such that $\alpha_i - 
\alpha_j \in \Phi$. By repeated application of Lemma \ref{lem1} (i)
it follows that $t_{\alpha_i} t_{\alpha_j} \preceq \gamma$. From $\alpha_i 
- \alpha_j \in \Phi$ we conclude that $\{\alpha_i, \alpha_j\}$ cannot be 
the simple system in the rank two induced subsystem of $\Phi$ spanned by 
$\alpha_i$ and $\alpha_j$. Since $<_\gamma$ is compatible with $\gamma$, 
Example \ref{ex:dihedral} shows that $t_{\alpha_i} >_\gamma t_{\alpha_j}$ 
which implies that $c$ is not rising.  
\qed

\begin{corollary}
If $W$ has type $A_{n-1}, B_n$ or $D_n$ then the natural edge labeling of 
$\nc_W (\gamma)$ is an EL-labeling under the choices of Coxeter element 
and total ordering on $T$ described in Examples \ref{ex:A} and \ref{ex:BD}. 
\qed
\label{cor:ABD}
\end{corollary}

\section{Proof of Theorem \ref{thm:main}}
\label{stein}

Let $W$ be any finite real reflection group of rank $\ell$ with set 
of reflections $T$, root system $\Phi$ and fixed choice of a positive 
system $\Phi^+$. Let $\{\sigma_1, \sigma_2,\dots,\sigma_\ell\}$ be a 
choice of simple system for $\Phi$ such that $\{\sigma_1,\dots,\sigma_r\}$ 
and $\{\sigma_{r+1},\dots,\sigma_\ell\}$ are orthonormal sets for some 
$r$ \cite[Section 3.17]{Hu} \cite{Ste}. 
It is proved in \cite{Ste} (under the additional assumption that $W$ is 
irreducible, which is acually not needed here) that $\Phi^+ = \{ \rho_1, 
\rho_2,\dots,\rho_{\ell h/2} \}$, where $\rho_i = t_{\sigma_1} 
t_{\sigma_2} \cdots t_{\sigma_{i-1}} (\sigma_i)$ for $1 \le i \le \ell 
h/2$ and the $\sigma_j$ are indexed cyclically modulo $\ell$. Consider 
the total ordering 
\begin{equation}
t_{\rho_1} < t_{\rho_2} < \cdots < t_{\rho_{\ell h / 2}} 
\label{steinberg}
\end{equation}
of $T$. The next statement follows from Theorem 5.4 in \cite{BW2}. 
\begin{theorem} {\rm (\cite{BW2})}
The total ordering {\rm (\ref{steinberg})} of $T$ is a reflection ordering 
for $W$ which is compatible with the Coxeter element $\gamma = t_{\sigma_1} 
t_{\sigma_2} \cdots t_{\sigma_\ell}$. 
\label{thm:weyl}
\end{theorem}

The previous theorem establishes the existence of a reflection ordering
which is compatible with some Coxeter element for any $W$. Combined with
Theorem \ref{thm:el} (ii) it gives a case-free proof of the statement in 
Theorem \ref{thm:main} in the case of Weyl groups. Using other results 
from \cite{BW2} we can give a different case-free proof that the ordering 
(\ref{steinberg}) yields an EL-shelling of $\nc_W (\gamma)$ for any $W$ 
as follows.

\begin{theorem}
If $T$ is totally ordered by {\rm (\ref{steinberg})} and $\gamma = 
t_{\sigma_1} t_{\sigma_2} \cdots t_{\sigma_\ell}$ then the natural edge 
labeling of $\nc_W (\gamma)$ with label set $T$ is an EL-labeling. 
\label{thm:any}
\end{theorem}
\begin{proof} 
We claim that there is at most one rising maximal chain with respect to 
the natural edge labeling $\lambda$ in any non-singleton interval $[u, v]$ 
in $\nc_W (\gamma)$. In view of part (i) of Theorem \ref{thm:el} it suffices 
to prove this claim. As in the proof of part (ii) of the same result, it 
suffices to treat the intervals of the form $[1, w]$. 

Let $(t_1, t_2,\dots,t_k)$ be the label of a rising maximal chain in 
$[1, w]$ and let $t_i = t_{\tau_i}$ for $1 \le i \le k$, where $\tau_i 
\in \Phi^+$. We will prove that $\{\tau_1, \tau_2,\dots,\tau_k\}$ is 
the set of simple roots of the subsystem $\Phi_w \subseteq \Phi$ 
corresponding to the parabolic subgroup $W_w$ (see Lemma \ref{lem5}) 
with respect to the positive system $\Phi_w \cap \Phi^+$. This clearly 
implies the claim. Let $\tau$ be the largest positive root in $\Phi_w$ 
with respect to (\ref{steinberg}). By part (i) of Lemma \ref{lem:basis} 
and Lemma \ref{lem2} the set $\{\tau_1, \tau_2,\dots,\tau_k\}$ is a 
basis of the real vector space spanned by $\Phi_w$. Hence there is a 
unique expression of the form  
\begin{equation}
\tau = a_1 \tau_1 + a_2 \tau_2 + \cdots + a_k \tau_k 
\label{eq}
\end{equation}
with $a_i \in \RR$ for all $i$. As in \cite[Lemma 3.9]{BW2} let $\mu 
(x) = -2 (\gamma - I)^{-1} x$, where $I$ is the identity map. For 
notational convenience we write $x \cdot y$ instead of $(x, y)$. From 
$\gamma = t_1 t_2 \cdots 
t_k$ and part (ii) of this lemma we have $\mu (\tau_i) \cdot \tau_j = 
0$ for $i < j$. Moreover \cite[Theorem 3.7]{BW2} implies that $\mu 
(\tau_j) \cdot \tau_i \le 0$ for $i<j$ and that $\mu (\tau_i) \cdot 
\tau_i = 1$ and $\mu (\tau_i) \cdot \tau \ge 0$ for all $i$. Taking the 
inner product of (\ref{eq}) with $\mu (\tau_i)$ we get
\begin{eqnarray*}
0 & \le & \mu(\tau_i) \cdot \tau \\
&=& a_1 \mu(\tau_i) \cdot \tau_1 +
\cdots + a_{i-1} \mu(\tau_i) \cdot \tau_{i-1} + a_i.
\end{eqnarray*}
It follows by induction that $a_i \ge 0$ for all $1 \le i \le k$. Thus 
$\tau$ is in the positive 
cone of $\{\tau_1, \tau_2,\dots,\tau_k\}$. Corollary 3.8 in \cite{BW2} 
gives $\tau_k = \tau$. By induction on the length of $w$ we may assume 
that $\{\tau_1, \tau_2,\dots,\tau_{k-1}\}$ is the simple system in 
$\Phi_{w t_k} \cap \, \Phi^+$. Theorem 5.1 in \cite{BW2} implies that 
$\{\tau_1, \tau_2,\dots,\tau_k\}$ is the simple system of $\Phi_w$, as 
desired. 
\end{proof}

\vspace{0.1 in}
\noindent
\emph{Proof of Theorem \ref{thm:main}.} It follows from Theorem
\ref{thm:any}.
\qed

\bigskip
The next corollary follows from Theorem \ref{thm:any} and standard 
facts on M\"obius functions of EL-shellable posets; see \cite[Section 
3.13]{Sta}. In view of Theorem \ref{thm:el} (ii) it applies to any 
Coxeter element and compatible reflection ordering in the case of 
Weyl groups.
\begin{corollary}
If $T$ is totally ordered by {\rm (\ref{steinberg})} then the M\"obius 
function on any interval $[u,v]$ in $\nc_W (\gamma)$ is equal to 
$(-1)^{l_T(v) - l_T(u)}$ times the number of falling maximal chains in 
$[u,v]$ with respect to the natural edge labeling of $\nc_W (\gamma)$. 
\qed
\label{cor:mobius}
\end{corollary}
In particular the M\"obius number of $\nc_W (\gamma)$ is equal to 
$(-1)^\ell$ times the number of falling maximal chains in $\nc_W (\gamma)$ 
with respect to this labeling. It can be deduced from the results of 
\cite[Section 8]{BW2} that this number of falling chains is equal to the 
number of positive clusters of the generalized associahedron of 
corresponding type.

\vspace{0.2 in}
\noindent
\emph{Acknowledgements}. 
The first two authors thank Jon McCammond, Alexandru Nica and Victor 
Reiner for their invitation to the AIM workshop \emph{Braid groups,
Clusters and Free Probability} in January 2005, in which some of the 
present work was done, the Institute Mittag-Leffler (Stockholm) and 
the Centre de Recerca Matem\`{a}tica (Barcelona), respectively, 
where this work was completed and 
Hugh Thomas for interesting discussions. The hospitality, financial 
support and excellent working conditions offered by the above 
mentioned institutions are gratefully acknowledged.

\end{document}